\documentclass[12pt]{amsart}
\newcommand{\sm}{\setminus}

\newcommand{\inv}{^{-1}}

\newcommand{\wt}{\widetilde}

\newcommand{\PP}{{\mathbb P}}
\newcommand{\R}{{\mathbb R}}
\newcommand{\C}{{\mathbb C}}

\newcommand{\half}{{\frac{1}{2}}}

\renewcommand{\phi}{\varphi}

\newcommand{\ocal}{\mathcal{O}}

\newcommand{\al}{\alpha}

\newcommand{\Ga}{\Gamma}

\newcommand{\La}{\Lambda}
\newcommand{\la}{\lambda}
\newcommand{\ep}{\varepsilon}
\newcommand{\de}{\delta}
\newcommand{\De}{\Delta}

\newtheorem{thm}{{Theorem}}[section]

\newtheorem{lem}[thm]{{Lemma}}
\newtheorem{prop}[thm]{{Proposition}}
\newenvironment{rmk}{\medskip\noindent{\it Remark:\/}}{\medskip}

\begin{document}

\title[Constructing low degree hyperbolic surfaces in $\PP^3$]{Constructing
low degree hyperbolic\\ surfaces in $\PP^3$}

\author{Bernard Shiffman}
\address{Department of Mathematics, Johns Hopkins University, Baltimore,
MD 21218, USA}
\email{shiffman@math.jhu.edu}

\author{Mikhail Zaidenberg}
\address{Universit{\'e}
Grenoble I, Institut Fourier, UMR 5582 CNRS-UJF, BP 74,
38402 St.\ Martin
d'H{\`e}res c{\'e}dex, France}
\email{zaidenbe@ujf-grenoble.fr}

\thanks{Research of the first author partially supported by NSF grants
\#DMS-9800479 and \#DMS-0100474.}



\begin{abstract}  We describe a new method of constructing
Kobayashi-hyperbolic
surfaces in complex
projective 3-space based on deforming surfaces with a ``hyperbolic
non-percolation'' property.  We use this method to show that general small
deformations of certain singular abelian surfaces of  degree 8 are
hyperbolic. We also show that a union of 15
planes in general position in projective 3-space admits hyperbolic
deformations.\end{abstract}
 
\maketitle

\section{Introduction}  A compact complex manifold is hyperbolic
in the sense of Kobayashi if  every holomorphic map  from the
complex line $\C$ to the manifold is constant, as is the case for
compact complex curves of genus $\ge 2$. In 1970, S. Kobayashi
conjectured that  generic hypersurfaces in $\PP^n$ of
sufficiently high degree are hyperbolic. Some progress has been
made towards this conjecture.  Demailly and El Goul \cite{DE},
and independently McQuillan \cite{MQ} (with a slightly bigger
degree estimate)  proved  that a very generic surface of
degree at least 21 in $\PP^3$ is hyperbolic in the sense of
Kobayashi. Previously, Clemens \cite{Cl} showed that very generic
hypersurfaces of degree $d\ge 2n-1$ contain no rational curves,
which is a necessary condition for hyperbolicity. Actually, for  $n\ge 3$,
this also holds for
$d= 2n-2$, and very generic hypersurfaces of degree $d\ge 2n-1$ contain
neither rational nor elliptic curves (Voisin \cite{Vo}; see also
\cite{CR, CLR, E, Pa, Xu1, Xu2}). Thus it is natural to suppose
that $2n-1$ is the minimal degree for Kobayashi's conjecture.

Many examples have been given of low degree hyperbolic projective
hypersurfaces (e.g., \cite{ShZa1,ShZa2} and the references therein).
The examples of hyperbolic surfaces in $\PP^3$ of lowest degree found
to date are of degree 8   and were discovered independently by  Duval
\cite{Du} and Fujimoto
\cite{Fu}.  (A similar example of degree 10 was previously
found by Shirosaki \cite{Shr}.)

In this paper, we provide a new approach to
constructing hyperbolic surfaces in $\PP^3$, giving another example of degree
8.  Although this example is considerably more complicated than the
Duval-Fujimoto example, we hope that our technique can be applied in the
future to construct examples of lower degree.  We
also show that  certain small deformations of 15 planes in
general position in $\PP^3$  are hyperbolic  surfaces.

Our technique involves showing that small deformations $X_t$ of
certain singular surfaces $X_0\subset \PP^3$ are hyperbolic.  The
surfaces $X_0$ that we deform, while not hyperbolic, satisfy a
``hyperbolic non-percolation" property. In particular, we
consider surfaces $X_0$ with ``double curve'' $\bar S$; i.e., the
singular locus of $X_0$ is a 1-dimensional subvariety $\bar S$,
and $X_0$ has two branches at general points of $\bar S$.  If
nearby surfaces in a linear pencil $\{X_t\}$ were not hyperbolic,
then we can find a sequence $t_n\to 0$ and a sequence of Brody
curves $f_{t_n}:\C\to X_{t_n}$ converging to a Brody curve
$f:\C\to X_0$.  Recall that a {\it Brody curve\/} in a Hermitian complex
manifold $M$ is a non-constant entire holomorphic curve $g:\C\to
M$ such that $\|g'(\zeta)\|$ is bounded above by $\|g'(0)\|$.
Brody \cite{Br} proved that a compact complex manifold is
hyperbolic iff it does not contain any Brody curves.

Our approach is to show using  Hurwitz's theorem that either
\begin{itemize}
\item $f(\C)\subset
\bar S\sm \{p_j\}$, where the $p_j$ are the multiple points of $\bar S$, or
\item $f(\C)\subset (X_0\sm \bar S)\cup D$, where $D$ is a finite
subset of $\bar S$.
\end{itemize}
We say
that  $X_0\sm \bar S$ {\it has the property of hyperbolic non-percolation
through\/} $D$ if there are no Brody curves $g:\C\to (X_0\sm \bar S)\cup D$.
Hence, if in addition,
\begin{itemize}
\item $\bar S\sm \{p_j\}$ is hyperbolic, and
\item $X_0\sm \bar S$ has
the property of hyperbolic non-percolation through $D$,
\end{itemize}
then small
deformations $X_t$ are hyperbolic. We illustrate our construction with two
examples.

\section{Deformation of 15 planes}

In 1989, the second author \cite{Za} showed that the complements of certain
smooth, irreducible small deformations of 5 lines in $\PP^2$ are complete
hyperbolic and hyperbolically embedded.  We begin by using our technique to
give, as a parallel example,  hyperbolic deformations of 15 planes in general
position in
$\PP^3$.

Let $L_j, \ j=1,\dots,15$, be linear functions on $\C^4$ defining
hyperplanes $$H_j:=\{z\in\PP^3:L_j(z)=0\}$$ in general position; i.e., any
4 of the $L_j$ are linearly independent, or equivalently, every point of
$\PP^3$ is contained in at most 3 of the $H_j$.   Let
$D=\{z\in\PP^3: Q(z)=0\}$ be a general quintic, and consider the linear pencil
of surfaces:
$$\textstyle G_t = \left\{\prod_{j=1}^{15}L_j + t
Q^3=0\right\}\subset\PP^3\;.$$

\begin{thm} The surface  $G_t$ is hyperbolic for sufficiently small $t\ne 0$.
\label{15planes} \end{thm}

\begin{proof} Suppose on the contrary that there exist $t_n\to 0$ such
that $G_{t_n}$ is not hyperbolic.  Then we can find a sequence of Brody
curves
$f_n:\C\to G_{t_n}$ with $\sup\|f_n'\|=\|f_n'(0)\|= 1$, where the norm is
computed with respect to  the Fubini-Study metric on $\PP^3$. Then we can
choose a subsequence, which we also denote by $\{f_n\}$, converging to a
Brody curve
$f:\C\to \bigcup_{j=1}^{15}H_j$. Assume without loss of generality that
$f(\C)\subset H_{15}$.

We first show that \begin{equation}\label{claim} f(\C)\subset \left(H_{15}\sm
\bigcup _{j=1}^{14} H_j \right)\cup D\;. \end{equation}

To verify (\ref{claim}), suppose  that $D$ does not pass through
any of the points $H_i\cap H_j\cap H_k$ ($i,j,k$ distinct), and
that on the contrary  $f(\zeta_0)\in H_k\sm D$, where $k\ne 15$.
Let $\De$ be a small disk about $\zeta_0$ so that $f(\bar \De)$
does not intersect $D$; i.e., $Q\circ f(\zeta)\ne 0$ for
$\zeta\in\bar\De$. Hence
$$\textstyle\prod_j L_j\circ f_n(\zeta)=-t_nQ^3\circ f_n(\zeta)\ne
0,\qquad \zeta\in\De,\ n\gg 0\;.$$  Since $L_k\circ
f(\zeta_0)=0$, it follows from Hurwitz's Theorem that $L_k\circ
f\equiv 0$; i.e., $f(\C)\subset H_k\cap H_{15}$.  If $j\not\in\{
k,15\}$ and $f(\C)$ passes through the point  $ p_{jk}:=H_j\cap
H_k\cap H_{15}\not\in D$, then by the above argument (replacing $k$
with $j$), we again conclude from Hurwitz's Theorem that
$f(\C)\subset
\{p_{jk}\}$, contradicting the fact that $f$ is  non-constant.  Hence
$$f(\C)\subset H_k\cap H_{15} \sm \{p_{jk}:1\le j\le 14, j\ne k\} \approx
\PP^1 \sm \{\mbox{13 points}\}\;,$$
which implies that $f$ is constant, a contradiction. Therefore (\ref{claim})
holds.

We assume that the curves $D\cap H_{j}$ are smooth
(or have at most 4 double points), so that the degree 5 curve
$D\cap H_{15}$ is hyperbolic and hence
$f(\C)\not\subset D$.  Then by the Cartan Second Main Theorem \cite{Ca} (see
also \cite[\S 3.B]{Ko}) applied to the map $f:\C \to
H_{15}\approx\PP^2$ and the 14 lines
$H_j\cap H_{15}$, we have
\begin{equation}\label{ineq1} (14-3)T_f(r) \le \sum_{j=1}^{14}N_2(H_j,r)
+O(\log T_f(r))\;.\end{equation} (Note that since Brody curves are of finite
order
$\le 2$, the inequality holds without exceptional intervals.) As we have
assumed that $D$ does not pass through any of the points $H_i\cap H_j\cap
H_k$, (\ref{claim}) implies that
\begin{equation}\label{ineq2}\sum_{j=1}^{14}N_2(H_j,r)  \le 2
\sum_{j=1}^{14}N_1(H_j,r)\le 2N(D,r)\;.\end{equation}
 Furthermore, by the Carlson-Griffiths First
Main Theorem \cite{CG} (see also \cite[\S 5.2]{NO}, \cite[\S
8.4]{Ko}) applied to the divisor
$D$, we have

\begin{equation}\label{ineq3}N(D,r)\le
5T_f(r)+O(1)\;.\end{equation}

Combining the inequalities (\ref{ineq1})--(\ref{ineq3}), we arrive at a
contradiction.

\end{proof}

\begin{rmk} In the second part of the proof of Theorem \ref{15planes}, we
showed that the complement of 14 general lines in $\PP^2$ has the property
of  hyperbolic non-percolation through (the intersection of these lines with)
a general quintic.  This should also hold for fewer lines; e.g., an open
problem is whether the complement of 5 general lines in
$\PP^2$  has the property
of  hyperbolic non-percolation through  a general
sextic curve.  This would imply that a general small deformation of 6
planes in general position in $\PP^3$ is a hyperbolic sextic
surface.\end{rmk}

\section{Deformation of a singular abelian surface}

We now use our hyperbolic non-percolation technique to
construct a new example of a hyperbolic surface of degree 8.  This time,
instead of deforming a reducible surface, we deform an irreducible surface
$X_0$ with self-intersections.

The surface $X_0$ is described in
\cite{LB} and is defined as follows.  Let $A$ be
a simple abelian surface with an ample line bundle  $L\to A$ of type
$(1,4)$. Recall that an abelian variety is said to
be simple if it does not contain any proper abelian subvariety.
See \cite[p.~47]{LB} for the definition of line bundles of type
$(d_1,d_2)$ on an abelian surface.

It follows
from the Riemann-Roch Theorem and the Kodaira vanishing theorem that
$h^0(A,L)=\chi(L)=\half L\cdot L =4$ and hence $\deg L:= L\cdot
L=8$ (see \cite[p. 289]{LB}).  We let
$X_0=\varphi_L(A)$, where $\varphi_L:A\to\PP^3$ is the rational map defined
by the linear system $|L|$.
We shall establish the following result:

\begin{thm}\label{thm8} General small deformations
of the surface $X_0\subset\PP^3$ are hyperbolic surfaces of
degree $8$.\end{thm}

We begin with the description of $X_0$.  By \cite[pp.~308--312]{LB}, the
surface $X_0$ is birational to $A$ and is given by
$$X_0=\{z\in\PP^3:Q=0\}\;,$$ where
\begin{eqnarray*}Q &=& {\lambda_{{1}}}^{2}\left
({z_{{0}}}^{4}{z_{{1}}}^{4}+{z_{{2}}}^{4}{z_{
{3}}}^{4}\right )+{\lambda_{{2}}}^{2}\left ({z_{{0}}}^{4}{z_{{2}}}^{4}
+{z_{{1}}}^{4}{z_{{3}}}^{4}\right )+{\lambda_{{3}}}^{2}\left ({z_{{0}}
}^{4}{z_{{3}}}^{4}+{z_{{1}}}^{4}{z_{{2}}}^{4}\right )\\
&&+2\,\lambda_{{1}} \lambda_{{2}}\left
({z_{{0}}}^{2}{z_{{1}}}^{2}+{z_{{2}}}^{2}{z_{{3}}}^ {2}\right
)\left (-{z_{{0}}}^{2}{z_{{2}}}^{2}+{z_{{1}}}^{2}{z_{{3}}}^{
2}\right )\\&&+2\,\lambda_{{1}}\lambda_{{3}}\left
({z_{{0}}}^{2}{z_{{1}}}^ {2}-{z_{{2}}}^{2}{z_{{3}}}^{2}\right
)\left ({z_{{0}}}^{2}{z_{{3}}}^{2
}-{z_{{1}}}^{2}{z_{{2}}}^{2}\right )\\&&
+2\,\lambda_{{2}}\lambda_{{3}} \left
({z_{{0}}}^{2}{z_{{2}}}^{2}+{z_{{1}}}^{2}{z_{{3}}}^{2}\right )
\left
({z_{{0}}}^{2}{z_{{3}}}^{2}+{z_{{1}}}^{2}{z_{{2}}}^{2}\right )+{
\lambda_{{0}}}^{2}{z_{{0}}}^{2}{z_{{1}}}^{2}{z_{{2}}}^{2}{z_{{3}}}^{2}\;.
\end{eqnarray*} In fact,  general choices of $(\la_0,\dots,\la_3)\in\C^4$
give singular abelian surfaces (see Remark 3.3 in
\cite[p.~301]{LB}).

Let $H_j$ denotes the coordinate plane $\{z_j=0\}$, and let 
$$p_j:=(\de^j_0:\de^j_1:\de^j_2:\de^j_3)\in X_0\qquad
\mbox{(so that\ }\ p_0=(1:0:0:0),\ \mbox{etc.)}$$ 
denote the vertices of the
coordinate tetrahedron $\{H_0,H_1,H_2,H_3\}$ ($0\le j\le 3$). The
singular locus of $X_0$ consists of 4 double curves $\bar
S_j:=X_0\cap H_j$, $j=0,1,2,3$. The equation for, say,  $\bar S_3\subseteq
H_3$ is
$$-\lambda_{{1}}{z_{{0}}}^{2}{z_{{1}}}^{2}+
\lambda_{{2}}{z_{{0}}}^{2}{z_{{2}}}^{2}+
\lambda_{{3}}{z_{{1}}}^{2}{z_{{2}}}^{2}=0\,.$$ It is
known \cite[p.~312]{LB} that $\bar S_3$ is an irreducible
rational curve with $3$ ordinary double points at
$\{p_0,p_1,p_2\}$.  Generic points of $\bar S_3$ are ordinary
double points of the surface $X_0$; i.e., $X_0$ is the union of
two transversal smooth surface germs at generic points of $\bar
S_3$. The set of points of $  \bar S_3$ which are not ordinary
double points of $X_0$ consists of the 3 double points
$\{p_0,p_1,p_2\}$ of $\bar S_3$ together with 12 smooth points of
$\bar S_3$, which are pinch points of $X_0$  (see
\cite[p.~312]{LB}). The same description applies to the other
double curves $\bar S_0,\bar S_1,\bar S_2$.

 We need to
know the structure of $X_0$ at the 4 vertices $\{p_j\}$. Let
$\bar S = \bigcup _{j=0}^3\bar S_j$.  The singular set of $\bar
S$ consists of the 4 points $\{p_j\}$. We shall show that each  $p_j$ is an
ordinary 6-fold singularity of $\bar S$; i.e., the germ of $\bar S$ at
$p_j$ consists of 6 smooth local curves $\{B_j^i\}_{1\le i\le 6}$ with
distinct tangents ($j=0,1,2,3$).  For example, $\bar S_1,\bar
S_2,\bar S_3$ pass through $p_0$, each contributing 2 local
components of the germ of $\bar S$ at $p_0$.
To describe the tangents to the
$B_0^i$, we write
$$\mu_1=\sqrt{\la_1}\,,\qquad  \mu_2=\sqrt{\la_2}\,,\qquad
\mu_3=\sqrt{-\la_3}\,,$$
and we use the
affine coordinates
$$x=z_1/z_0, \ y=z_2/z_0,\ z=z_3/z_0$$ about $p_0=(0,0,0)$. Then the 6
tangent lines are
$$\ell_1^\pm = \{x=0\,,\ \mu_2y=\pm\mu_3 z\},\ \
\ell_2^\pm = \{y=0\,,\ \mu_1x=\pm\mu_3 z\},\ \
\ell_3^\pm = \{z=0\,,\ \mu_1x=\pm\mu_2y\}\,,$$
where the 2 lines $\ell_j^+,\ell_j^-$ are tangent to $\bar S_j$ at $p_0$
($j=1,2,3$).

\begin{lem}\label{br} Let $0\le j\le 3$.  The germ of $X_0$ at $p_j$
consists of 4 smooth surface germs $Y_j^1,\dots, Y_j^4$.  Each of these
surfaces contains exactly 3 of the 6 components $B_j^1,\dots B_j^6$ of
the germ of $\bar
S$ at ${p_j}$, and each of these components  $B_j^k$ is
contained in exactly 2 of
the $Y_j^i$.  The intersection of any 3 of the $Y_j^i$ is the
germ of the point $p_j$.\end{lem}

\begin{proof} Clearly,  the orbit of $p_0$ under
Aut$(X_0)$ consists of the 4 vertices $\{p_j\}$.
Thus it suffices to consider $j=0$. As before, fix  the  affine
coordinates
$$x=z_1/z_0, \ y=z_2/z_0,\ z=z_3/z_0$$ about $p_0$.  Then
\begin{eqnarray*}Q&=& {\lambda_{{1}}}^{2}\left ({x}^{4}
+{y}^{4}{z}^{4}\right
)+{\lambda_{{2} }}^{2}\left ({y}^{4}+{x}^{4}{z}^{4}\right
)+{\lambda_{{3}}}^{2}\left ( {z}^{4}+{x}^{4}{y}^{4}\right
)\\&& +2\,\lambda_{{1}}\lambda_{{2}}\left ({x} ^{2}+{y}^{2}{z}^{2}\right
)\left (-{y}^{2}+{x}^{2}{z}^{2}\right ) +2\,
\lambda_{{1}}\lambda_{{3}}\left ({x}^{2}-{y}^{2}{z}^{2}\right )\left (
{z}^{2}-{x}^{2}{y}^{2}\right )\\&&
+2\,\lambda_{{2}}\lambda_{{3}}\left ({y}
^{2}+{x}^{2}{z}^{2}\right )\left ({z}^{2}+{x}^{2}{y}^{2}\right )+{
\lambda_{{0}}}^{2}{x}^{2}{y}^{2}{z}^{2}\\&=& az^4 +bz^2+c\,,
\end{eqnarray*} where
\begin{eqnarray*}a&=&{\lambda_{{2}}}^{2}{x}^{4}+2\,\lambda_{{2}}\lambda_{{1}}{x}^{2}{y}^{2}
+{\lambda_{{1}}}^{2}{y}^{4}+2\,\lambda_{{3}}\lambda_{{2}}{x}^{2}-2\,
\lambda_{{3}}\lambda_{{1}}{y}^{2}+{\lambda_{{3}}}^{2}\,,\\
b&=&2\,\lambda_{{3}}\lambda_{{2}}{x}^{4}{y}^{2}+2\,\lambda_{{3}}\lambda_{{
1}}{x}^{2}{y}^{4}+2\,\lambda_{{2}}\lambda_{{1}}{x}^{4}-2\,\lambda_{{2}
}\lambda_{{1}}{y}^{4}+{\lambda_{{0}}}^{2}{x}^{2}{y}^{2}+2\,\lambda_{{3
}}\lambda_{{2}}{y}^{2}+2\,\lambda_{{3}}\lambda_{{1}}{x}^{2}
\,,\\
c&=&{\lambda_{{3}}}^{2}{x}^{4}{y}^{4}+2\,\lambda_{{3}}\lambda_{{2}}{x}^{2}
{y}^{4}-2\,\lambda_{{3}}\lambda_{{1}}{x}^{4}{y}^{2}+{\lambda_{{2}}}^{2
}{y}^{4}-2\,\lambda_{{2}}\lambda_{{1}}{x}^{2}{y}^{2}+{\lambda_{{1}}}^{
2}{x}^{4}\\
&=&\left (\lambda_{{3}}{x}^{2}{y}^{2}-\lambda_{{1}}{x}^{2}+\lambda_{{2}}{
y}^{2}\right )^{2}
\,.\end{eqnarray*}
We compute the discriminant
$$\de :=b^2-4ac=x^2y^2\left(16\la_1\la_2{\la_3}^2 + \cdots\right)\;.$$
Hence in  the local ring ${}_3\ocal_0=\C[[x,y,z]]$, we have the factorization
$Q=Q^+Q^-$, where
\begin{eqnarray*}Q^\pm\ =\ \sqrt{a}\,z^2 +\frac{b\pm\sqrt{\de}}{2\sqrt{a}
}&=&\la_3 z^2+\la_1 x^2 +\la_2y^2 \pm 2\sqrt{\la_1\la_2}\,xy+R_3(Q^\pm)\\
&=& (\mu_1x \pm \mu_2 y)^2 -(\mu_3z)^2+R_3(Q^\pm)\,.
\end{eqnarray*} Here, for any $f\in{}_3\ocal_0$, we let
$R_n(f)=\sum_{k=n}^\infty f_k$, where $f_k$ denotes the homogeneous terms of
order $k$ in the Taylor expansion of $f$.

It follows that the tangent cone of $X_0$ at $p_0$ is the union of the 4
planes (in general position):
\begin{eqnarray*}
P^{++} &=& \{\mu_1x+\mu_2y+\mu_3z=0\}\,,\\
P^{+-} &=& \{\mu_1x+\mu_2y-\mu_3z=0\}\,,
\\P^{-+} &=& \{\mu_1x-\mu_2y+\mu_3z=0\}\,,\\
P^{--} &=& \{\mu_1x-\mu_2y-\mu_3z=0\}\,.
\end{eqnarray*}
The planes $P^{++}$ and $P^{+-}$  are tangent to $\{Q^+=0\}$, while
$P^{-+}$ and $P^{--}$  are tangent to $\{Q^-=0\}$.

To show that $Q^+$ and $Q^-$ are reducible in ${}_3\ocal_0$, we now expand
$Q=a'x^4+b'x^2+c'$, where
\begin{eqnarray*}a'&=&{\lambda_{{3}}}^{2}{y}^{4}+2\,\lambda_{{3}}\lambda_{{2}}{z}^{2}{y}^{2}
+{\lambda_{{2}}}^{2}{z}^{4}-2\,\lambda_{{3}}\lambda_{{1}}{y}^{2}+2\,
\lambda_{{2}}\lambda_{{1}}{z}^{2}+{\lambda_{{1}}}^{2}
\,,\\
b'&=&2\,\lambda_{{3}}\lambda_{{1}}{z}^{2}{y}^{4}+2\,\lambda_{{2}}\lambda_{{
1}}{z}^{4}{y}^{2}+2\,\lambda_{{3}}\lambda_{{2}}{z}^{4}+2\,\lambda_{{3}
}\lambda_{{2}}{y}^{4}+{\lambda_{{0}}}^{2}{z}^{2}{y}^{2}+2\,\lambda_{{3
}}\lambda_{{1}}{z}^{2}-2\,\lambda_{{2}}\lambda_{{1}}{y}^{2}
\,,\\
c'&=&{\lambda_{{1}}}^{2}{z}^{4}{y}^{4}-2\,\lambda_{{3}}\lambda_{{1}}{z}^{4}
{y}^{2}-2\,\lambda_{{2}}\lambda_{{1}}{z}^{2}{y}^{4}+{\lambda_{{3}}}^{2
}{z}^{4}+2\,\lambda_{{3}}\lambda_{{2}}{z}^{2}{y}^{2}+{\lambda_{{2}}}^{
2}{y}^{4}
\\
&=&\left (-\lambda_{{1}}{z}^{2}{y}^{2}+\lambda_{{3}}{z}^{2}+\lambda_{{2}}
{y}^{2}\right )^{2}
\,.\end{eqnarray*}
Again we compute the discriminant
$$\de'=b'{}^2-4a'c'=y^2z^2\left(-16{\la_1}^2\la_2{\la_3}+
\cdots\right)\;.$$ Hence we have the factorization
$Q=Q'{}^+Q'{}^-$, where
\begin{eqnarray*}Q'{}^\pm\ =\ \sqrt{a'}\,x^2
+\frac{b'\pm\sqrt{\de'}}{2\sqrt{a'} }&=& (\mu_1x )^2
-(\mu_2y\pm\mu_3z)^2+R_3(Q'{}^\pm)\,.
\end{eqnarray*}
This time the planes $P^{++}$ and $P^{--}$  are tangent to $\{Q'{}^+=0\}$,
while
$P^{+-}$ and $P^{-+}$  are tangent to $\{Q'{}^-=0\}$.

Since ${}_3\ocal_0$ is a unique factorization domain, $Q^+$ and $Q'{}^+$ must
have a common factor $$Q^{++}=\mu_1x+\mu_2y+\mu_3z +R_2(Q^{++}) $$ with zero
set tangent to the plane
$P^{++}$.  By considering all such possible pairs, we see that the germ
of
$X_0$ at $p_0$ consists of 4 smooth surfaces, each tangent to one of the
planes $P^{++},P^{+-},P^{-+},P^{--}$.
One easily checks that each of the lines $\ell_j^\pm$ is the intersection of
exactly 2 of these planes, and each plane contains exactly 3 of the lines.
The conclusion of the lemma immediately follows.
\end{proof}

We recall that a simple complex torus does not contain
rational or elliptic curves. Moreover, the following holds.

\begin{prop}\label{A} Let  $f:\C\to T$ be a Brody curve in a simple
complex 2-dimensional torus $T=\C^2/\La$. Then for any compact
complex curve $S$ in $T$, the intersection $f(\C)\cap S$ is
infinite.
\end{prop}

\begin{proof}
Assume without loss of generality that $S$ is irreducible. The
lift $\tilde f: \C\to \C^2$ of $f$ is also a Brody curve and
hence is given by degree $1$ polynomials, as observed in
\cite{Gr}. By a translation of coordinates, we may suppose that
$\tilde f$ is linear and hence $f(\C)$ is a subgroup of $T$. We
first note that $S$ is not contained in any translate of
$f(\C)$.  For if on the contrary $S':=S+v\subset f(\C)$, then
$f\inv(S')$, being analytic and of positive dimension, must be
all of $\C$; i.e., $S'=f(\C)$. Hence $S'$, being a  compact
complex curve and a subgroup of $T$, must be a complex
1-dimensional subtorus of $T$, contradicting the assumption that
$T$ is simple.

Let $Y$ denote the closure of $f(\C)$ in the metric topology on
$T$.  As $T$ is a simple Lie group and $Y\subseteq T$ is a closed
subgroup, we  conclude that  either $Y=T$ or $Y$ is a real
3-subtorus of $T$. Choose new coordinates $\{z_1=x_1+iy_1,
z_2=x_2+iy_2\}$ in $\C^2$ so that $f(\C)$ is the image of the
axis $\{z_1=0\}$ via the projection $\tau:\C^2\to T$, and in the
latter case, $Y$ is the image of $\{y_1=0\}$.

\medskip

\noindent {\it Claim: \ $S\cap Y$ is
nonempty.}

\medskip \noindent {\it Proof of the claim:\/} We need to consider
only the case where $Y$ is a real 3-torus. Notice that the
universal cover of $T\sm Y$ can be identified with a strip in
$\C^2=\R^4$ between two parallel hyperplanes $y_1=0$ and
$y_1=a>0$, and that $y_1$ generates a well-defined bounded
harmonic function on $T\sm Y$. If $S\cap Y=\emptyset$, then
$y_1|_S=$const, whence $z_1|_S=$const, so $S$ is contained in a
translate of $f(\C)$, which is impossible as noted above. This
completes the proof of the claim.
\medskip

Now let
$\Delta_\ep^2$ be a  coordinate bidisk centered about a point
$s=(s_1,s_2)\in Y\cap S$.  We may assume by our choice of coordinates in $\C^2$
that $\Im s_1=0$ and
$$\Delta_\ep^2=\De'\times \De''=\{(w_1,w_2)\in\C^2: |w_j-s_j|<\ep,\;
j=1,2\}\buildrel {\tau}\over \hookrightarrow T\;.$$ We also
have
$f(\C)\cap
\Delta_{\ep}^2=E\times \De''$, where $E$ is a dense
subset of the disk
$\Delta'$ if $Y=T$, or is a dense subset of the real interval
$I_{\ep}:=\De'\cap\R$ if $Y$ is a real 3-torus.

We observe  that $S\not\supset \{s_1\}\times
\De''$.  Indeed, if on the contrary  $S$ contains the disk $\{s_1\}\times
\De''$, then a translate $S'$ of $S$ contains a disk $\{s_1'\}\times
\De''\subset f(\C)$, so the analytic set
$f\inv(S')$ cannot be 0-dimensional. Hence $f\inv(S')=\C$, or equivalently
$f(\C)\subset S'$, and therefore $Y=\overline{f(\C)}\subset S'$, a
contradiction.

Let $\rho_1:\Delta_\ep^2\to
\Delta'$ denote the projection to the
$z_1$-axis.
Since  $S\cap (\{s_1\}\times \De'')$ must be 0-dimensional, it follows that
$\rho_1(S\cap
\Delta_{\ep}^2)$ contains a neighborhood of $s_1$.  Since $s_1$ is a cluster
point of $E$, the set $E\cap\rho_1(S\cap
\Delta_{\ep}^2)$ must be infinite.  Since
$$E\cap\rho_1(S\cap
\Delta_{\ep}^2) = \rho_1((E\times \De'')\cap S)= \rho_1(f(\C)\cap S\cap
\De_\ep^2)\;,$$ it follows that $f(\C)\cap S$ is also infinite.
\end{proof}

\begin{rmk}\label{C} (i) Proposition \ref{A} can be rephrased as follows:
For any  compact  complex curve $S$ in a simple abelian surface
$T$ and for any divisor $D$ on $S$, the complement $T\backslash
S$ has the property of hyperbolic non-percolation through $D$.

(ii)  Note that if  $S$ were a rational or elliptic curve in a
simple complex torus $T$, then we would have a Brody curve
$f:\C\to S\subset T$, and by the first paragraph of the proof of
Proposition \ref{A}, $S$ must  contain a translate of a subgroup
of $T$, contradicting the assumption that $T$ is simple.\end{rmk}

\medskip

\noindent{\sc Proof of Theorem \ref{thm8}.} Let
$X_{\infty}=\{z\in\PP^3:F(z)=0\}$ be a general octic surface and consider the
linear pencil of surfaces $$X_t=\{z\in\PP^3:Q(z)+tF(z)=0\}\;.$$
Suppose on the contrary that $X_{t_n}$  is not hyperbolic for some sequence
$t_n\to 0$.  Then as in the proof of Theorem \ref{15planes}, after
passing to a subsequence of $\{t_n\}$, we can find a sequence of Brody curves
$f_n:\C\to X_{t_n}$ converging to  a Brody curve
$f:\C\to X_0$.

\medskip \noindent {\it Claim:\ \ $f(\C)\subset (X_0\backslash \bar
S)\cup (\bar S\cap X_{\infty})\cup \Ga$, where
$\Ga$ is the set of 48 pinch points of $X_0$.}

\medskip \noindent {\it Proof of the claim:\/}  Suppose on the contrary
that
$f(\zeta_0)=x_0\in \bar S\sm (X_\infty\cup \Ga)$. We first consider the case
where
$x_0\not\in
\{p_0,p_1,p_2,p_3\}$, so that $x_0$ is an ordinary double point of $X_0$.
Choose a small neighborhood
$U\subset \PP^3\sm X_\infty$ of $x_0$ such that we have a factorization
$Q|_U=Q'Q''$, where
$Q',\ Q''$ are holomorphic on $U$ and vanish at $x_0$; hence $X_0\cap U$
consists of two components $X'=\{Q'=0\}$, $X''=\{Q''=0\}$.
Let $\De$ be a small disk about $\zeta_0$ such that
$f(\bar\De)\subset U$.  Then (by the same argument as in the
proof of Theorem~\ref{15planes}) for $n$ sufficiently large, $f_{n}(\De)$
does not meet
$X_0$ and hence
$Q\circ f_n(\zeta)\ne 0$ for $\zeta\in \De$. Since
$f(\zeta_0)=x_0\in  X'\cap X''$,  it follows by Hurwitz's theorem applied to
$Q'\circ f_n$ and to $Q''\circ f_n$ that
$f(\De)\subset X'\cap X'' =\bar S\cap U$, and thus $f(\C)\subset \bar S$.
We shall complete this case below.

We now turn to the case $x_0=p_j$ ($j=0,1,2,3$).  By  virtue of
Lemma \ref{br}, this time  $ X_0\cap U$  consists of 4 components, and
we conclude as above that $f(\C)$ is contained in the intersection of these 4
components. But  (by  Lemma \ref{br}) this
intersection is the point $p_j$, and hence $f$ is constant, a
contradiction.

Returning to the first case, we can now conclude that
$f:\C\to\bar S\sm \{p_0,p_1,p_2,p_3\}$. However, the variety
$\bar S\sm \{p_0,p_1,p_2,p_3\}$ consists of 4 components $\bar
S_j\sm \{p_i:i\ne j\}$, $0\le j\le 3$, each a $\PP^1$ with 6
points (corresponding to the 3 double points of $S_j$) punctured
out.  Again we conclude that $f$ is constant, and this
contradiction completes the proof of the claim. 

\medskip
As the morphism  $\varphi_L: A \to X_0$ is birational and
proper, it provides a normalization of $X_0$. Let $\tilde f:\C\to
A$ be the lift of $f$ to the normalization, and let
$S:=\varphi_L\inv(\bar S)\subset A$.   The image $\tilde f(\C)$ meets
$S$ inside the finite set $D:=S\cap\varphi_L\inv(X_{\infty}\cup \Ga)$.
Although $f$ is a Brody curve, $\tilde f$ is not {\it a priori\/}
Brody, but we can construct a Brody curve from $\tilde f$ as
follows. Let $\De_n$ denote the disk of radius $n$ about the
origin.   By Brody's reparametrization lemma \cite{Br}, we  can  find a
sequence of holomorphic maps
$$g_n=\tilde f\circ \rho_n\circ \al_n:\De_n\to
\wt f(\C)\subset (A\sm S)\cup D\,,$$ with  $||{ g}_n'(0)||=c>0$
and $||{ g}_n'(\zeta)||\le \frac{n^2 c}{n^2-|\zeta|^2}$, where $\al_n$ is
a suitably chosen automorphism of $\De_n$ and $\rho_n(\zeta)=r_n \zeta$,
$r_n >0$ (see
\cite[(3.6.2), (3.6.4)]{Ko}, \cite[(1.6.6)]{NO}). 
After passing to a subsequence, we may suppose
that $ g_n$ converges as $n\to\infty$ to a Brody curve $ g :
\C\to A$.

By  Proposition \ref{A}, $g(\C)\cap
S$ is infinite.  Hence, there is a point $\zeta_0\in\C$ such that
$g(\zeta_0)\in S\sm D$. We choose a small disk $\De$ about $\zeta_0$ such
that $g_n(\De)$ does not meet $D$ and hence $g_n(\De)\subset A\sm S$ for $n\gg
0$. We then conclude as before by Hurwitz's theorem that
$g(\C)\subset S$.  But since $A$ is a simple abelian variety, $S$ cannot be
rational or elliptic.  Thus $g$ is constant, a contradiction.

Therefore,
$X_t$ is hyperbolic for $t$ sufficiently small.
\qed

\begin{rmk} Note that the only condition on $X_\infty$ is that
it does not contain any of the $p_j$. \end{rmk}

\noindent {\it Acknowledgment.\/} We are grateful to Gerd Dethloff for
useful comments. We also thank Herb Clemens for informing us of
the preprint \cite{ZR} as well as his own recent related
results (unpublished).

\end{document}